\documentclass[english,reprint,superscriptaddress,secnumarabic,nofootinbib,nobibnotes]{revtex4-1}
\usepackage[T1]{fontenc}
\usepackage[latin9]{inputenc}
\usepackage{geometry}
\usepackage{xcolor}
\usepackage{textcomp}
\usepackage{amsmath}
\PassOptionsToPackage{normalem}{ulem}
\usepackage{ulem}

\makeatletter

\providecommand{\tabularnewline}{\\}
\providecolor{lyxadded}{rgb}{0,0,1}
\providecolor{lyxdeleted}{rgb}{1,0,0}

\DeclareRobustCommand{\lyxsout}[1]{\ifx\\#1\else\sout{#1}\fi}

\setcitestyle{authoryear,round}

\usepackage{amsthm}

\providecommand{\tabularnewline}{\\}

\def\frontmatter@abstractheading{}

\renewcommand{\p@subsection}{}
\renewcommand{\p@subsubsection}{}

\usepackage{babel}

\makeatother

\usepackage{babel}
\begin{document}

\title{On Universality of Classical Probability with Contextually Labeled
Random Variables}

\author{Ehtibar N. Dzhafarov}

\thanks{Corresponding author: Ehtibar Dzhafarov, Purdue University, Department
of Psychological Sciences, 703 Third Street West Lafayette, IN 47907,
USA. email: ehtibar@purdue.edu.}

\affiliation{Purdue University}

\author{Maria Kon}

\affiliation{Purdue University}
\begin{abstract}
One can often encounter claims that classical (Kolmogorovian) probability
theory cannot handle, or even is contradicted by, certain empirical
findings or substantive theories. This note joins several previous
attempts to explain that these claims are unjustified, illustrating
this on the issues of (non)existence of joint distributions, probabilities
of ordered events, and additivity of probabilities. The specific focus
of this note is on showing that the mistakes underlying these claims
can be precluded by labeling all random variables involved contextually.
Moreover, contextual labeling also enables a valuable additional way
of analyzing probabilistic aspects of empirical situations: determining
whether the random variables involved form a contextual system, in
the sense generalized from quantum mechanics. Thus, to the extent
the Wang-Busemeyer QQ equality for the question order effect holds,
the system describing them is noncontextual. The double-slit experiment
and its behavioral analogues also turn out to form a noncontextual
system, having the same probabilistic format (cyclic system of rank
4) as the one describing spins of two entangled electrons. 
\end{abstract}
\maketitle
In the literature on foundations of\textcolor{black}{{} quantum physics
(Accardi, 1982; Feynman, 1951; Feynman, Leighton, \& Sands, 1975;
Khrennikov, 2009b)} and, more recently\textcolor{black}{, psychology
(Aerts, 2009, 2014; }Broekaert, Basieva, Blasiak, \& Pothos\textcolor{black}{,
2017; Busemeyer \& Bruza, 2012; Moreira \& Wichert, 2016; Pothos \&
Busemeyer, 2013)}, one can encounter statements that classical (Kolmogorovian)
probability theory does not have adequate conceptual means to handle
(sometimes, even, is contradicted by) this or that empirical fact. 

Three of the \textcolor{black}{most} widespread assertions of this
kind are as follows: 
\begin{description}
\item [{Statement$\:$1}] Classical probability requires that certain (e.g.,
Bell-type) inequalities hold for certain systems of random variables,
but we know from quantum mechanics and from behavioral experiments
that they may be violated. 
\item [{Statement$\:$2}] In classical probability, the joint occurrence
of two events is commutative, but we know from quantum mechanics and
from behavioral experiments that the order of two events generally
matters for their joint probability.
\item [{Statement$\:$3}] Classical probability is additive (equivalently,
obeys the law of total probability), but we know from quantum mechanics
and from behavioral experiments that this additivity (the law of total
probability) can be violated.
\end{description}
This note has three objectives: (1) to show that the three statements
above are based on misidentification of the random variables involved,
due to ignoring their inherently contextual labeling; (2) to show
that contextual labeling is a principled way to ``automatically''
ensure correct applicability of classical probability theory to an
empirical situation; and (3) to demonstrate how the use of contextual
labeling enables so-called contextuality analysis of systems of random
variables, a relatively new form of probabilistic analysis of considerable
interest in empirical applications. Contextual labeling of random
variables is the departing principle of Khrennikov's V\"axj\"o Model
(Khrennikov, 2009a) and of the Contextuality-by-Default theory (Dzhafarov,
2017; Dzhafarov, Cervantes, \& Kujala, 2017; Dzhafarov \& Kujala,
2014a, 2015, 2016a, 2017a, 2017b; Dzhafarov, Kujala, \& Cervantes,
2016). 

Let us preamble this discussion by stating our view of classical probability
theory (CPT), one that we are not prepared to defend in complete generality,
confining ourselves instead to merely illustrating it on the three
statements above. This view is that CPT, on a par with classical logic
and set theory, is a universal abstract mathematical theory. As an
abstract mathematical theory, it does not make empirically testable
predictions, because of which it cannot be contradicted by any empirical
observation. As a universal theory, for any empirical situation, it
has conceptual means to adequately describe anything that can be qualified
as this situation's probabilistic features (in the frequentist sense).
Moreover, as a conceptual tool, in the same way as classical logic
and set theory, it is indispensable and irreplaceable in dealing with
probabilistic problems: at the end, the results of any non-classical
probabilistic analysis have to be formulated in terms of classical
(frequentist) probabilities, distributions, and random variables.
However, when applied to an empirical situation, CPT can (even must)
be complemented by special-purpose computations identifying some of
the random variables, distributions, and probabilities in this particular
situation. To give a very simple example, CPT provides methods for
deriving probabilities of events defined on the outcomes of rolling
a die from a distribution of these outcomes, but it cannot predict
this distribution. A special theory is needed to know, e.g., that
if a die is manufactured in a particular way, then the distribution
of its outcomes is uniform. We view quantum probability as such a
special-purpose theory complementing classical probability. This mathematical
formalism is indispensable in quantum mechanics and has significant
achievements to its credit in psychology (e.g., Wang \& Busemeyer,
2013). It can be formalized and presented as an abstract calculus
alternative to or even generalizing the calculus of CPT, in the same
way one can formalize a paraconsistent logic as a generalization of
classical logic. However, just as one cannot replace classical logic
with paraconsistent logic in analyzing anything, including the very
paraconsistent logic itself, one cannot dispense with classical probability
when discussing and analyzing quantum probability computations and
relating them to data. 

This view is not entirely new. Ballentine (1986) defended a similar
position in essentially the same way we are doing here. The difference
is in that instead of using random variables, Ballentine confined
himself to a more limited language of events, and he used conditionalization
in place of the more general contextualization (Dzhafarov \& Kujala,
2014b; we discuss conditionalization in Section \ref{sec:On-Statement-3}
below). Khrennikov (2009a), in describing his V\"axj\"o contextual
model uses Ballentine's conditional-probability notation, but emphasizes
that these are not conditional probabilities of CPT. Rather he calls
them ``contextual probabilities,'' and explains that ``contextual
probability {[}...{]} is not probability that an event, say $B$,
occurs under the condition that another event, say $C$, has occurred.
The contextual probability is probability to get the result $a=\alpha$
under the complex of physical conditions $C$'' (\textcolor{black}{Khrennikov,
2009}a, p. 50). This seems to be the same as the contextual labeling
used in the Contextuality-by-Default theory. A very clear presentation
of a position that is close to ours can be found in the arguments
presented in an internet discussion by Tim Maudlin (2013). 

The purpose of this paper is to achieve conceptual clarity in understanding
CPT, not to criticize specific authors or papers. The latter is an
ungrateful task, as most authors' positions are not entirely consistent,
are subject to (re)interpretations, and evolve over time. We cite
specific papers and occasionally provide quotes only to demonstrate
that a reasonable reader may interpret the positions they entail in
the spirit of the Statements 1-3 above. Thus, Richard Feynman is often
cited as arguing that classical probability is not compatible with
q\textcolor{black}{uantum mechanics (Accardi, 1982; Costantini, 1993;
Khrennikov 2009b)}. This interpretation is supported by Feynman's
speaking of ``the discovery that in nature the laws of combining
probabilities were not those of the classical probability theory of
Laplace\textquotedblright{} (Feynman, 1951, p. 533). However, one
can also find statements in Feynman's writings that make his point
of view less than unequivocal. Thus, we read in the same paper and
on the same page that ``the concept of probability is not altered
in quantum mechanics. When I say the probability of a certain outcome
of an experiment is $p$ {[}...{]} no departure from the concept used
in classical statistics is required. What is changed, and changed
radically, is the method of calculating probabilities'' (\emph{ibid}).
This quote is consistent with treating quantum formalisms as special-purpose
computations embedded in CPT. We will return to Feynman when discussing
the double-slit experiment in Section \ref{sec:On-Statement-3}.

\section{\label{sec:On-Statement-1}On Statement 1}
\begin{quotation}
``Classical probability requires that certain (e.g., Bell-type) inequalities
hold for certain sets of random variables, but we know from quantum
mechanics and from behavioral experiments that they may be violated.'' 
\end{quotation}
\textcolor{black}{This view is commonly held in both physics and psychology
(Aerts, 2009; Aerts \& Sozzo, 2011; Bruza, Kitto, Nelson, \& McEvoy,
2009; Busemeyer \& Bruza, 2012; Filipp \& Svozil, 2005; Khrennikov,
2009b; Yearsley \& Pothos, 2014). In particular, among those applying
quantum probability to behavior and also treating quantum probability
theory as an alternative to CPT, there are claims that Bell-type inequalities
are violated in experiments involving combinations of concepts (Aerts
\& Sozzo, 2011; Busemeyer \& Bruza, 2012) and memory (Bruza, Kitto,
Nelson, \& McEvoy, 2009). }

We will not recapitulate all the arguments related to this issue,
as they have been presented in many previous publications (Dzhafarov,
Cervantes, \& Kujala, 2017; Dzhafarov \& Kujala, 2014a, 2014b, 2016a,
2017a, 2017b; Dzhafarov, Kujala, \& Larsson, 2015). We will use just
one familiar example. Let $R_{1},R_{2},R_{3},R_{4}$ denote a set
of binary ($+1/-1$) random variables with known distributions of
$\left(R_{1},R_{2}\right)$, $\left(R_{2},R_{3}\right)$, $\left(R_{3},R_{4}\right)$,
and $\left(R_{4},R_{1}\right)$. The necessary and sufficient condition
for the existence of such a quadruple of random variables is given
by the CHSH/Fine inequality (Bell, 1964; Clauser, Horne, Shimony,
\& Holt, 1969; Fine, 1982):
\begin{equation}
\max_{j=1,\ldots,4}\left|\sum_{i=1}^{4}\left\langle R_{i}R_{i\oplus1}\right\rangle -2\left\langle R_{j}R_{j\oplus1}\right\rangle \right|\leq2,\label{eq: CHSH}
\end{equation}
where $\oplus1$ is cyclic shift $1\rightarrow2\rightarrow3\rightarrow4\rightarrow1$,
and $\left\langle \cdot\right\rangle $ is expectation. One can easily
construct examples of distributions of $\left(R_{i},R_{i\oplus1}\right)$
for which this inequality is violated, indicating that such $R_{1},R_{2},R_{3},R_{4}$
do not exist (essentially by the same logic as in determining that
there are no four numbers $a,b,c,d$ with $a=b$, $b=c$, $c=d$,
and $d=a+1$).

The problem arises when we are being told that the existence of such
$R_{1},R_{2},R_{3},R_{4}$ is predicted by quantum theory and corroborated
by experiment. If we believe this, violations of ($\ref{eq: CHSH}$)
should indeed mean that CPT is inadequate, if not internally contradictory.
We should not, however, believe this. $R_{1},R_{2},R_{3},R_{4}$ in
($\ref{eq: CHSH}$) are random variables in the CPT sense; they are
not within the language of quantum theory. To decide what classical
random variables should describe outcomes of what quantum measurements,
one needs to go outside this theory. The general rule is that a random
variable is identified by what is being measured and how it is being
measured. The latter includes all conditions under which the measurement
is made, in particular, all other measurements performed together
with the given one. In our example, the measurements are indicated
by star symbols in the following matrix:
\begin{center}
\begin{tabular}{|c|c|c|c|c|}
\hline 
$\star$ & $\star$ &  &  & $c_{1}$\tabularnewline
\hline 
 & $\star$ & $\star$ &  & $c_{2}$\tabularnewline
\hline 
 &  & $\star$ & $\star$ & $c_{3}$\tabularnewline
\hline 
$\star$ &  &  & $\star$ & $c_{4}$\tabularnewline
\hline 
$q_{1}$ & $q_{2}$ & $q_{3}$ & $q_{4}$ & \multicolumn{1}{c}{}\tabularnewline
\cline{1-4} 
\end{tabular}$\:$.
\par\end{center}

\noindent The row labels $c_{1},\ldots,c_{4}$ are called contexts,
and here they are defined by which two quantities are being measured
together: in $c_{1}$ it is $q_{1}$ and $q_{2}$, in $c_{2}$ it
is $q_{2}$ and $q_{3}$, etc. In behavioral science the quantities
$q_{1},\ldots,q_{4}$ can be, e.g., four Yes-No questions posed to
a large number of people divided into four groups: in the group $c_{1}$
each person is asked $q_{1}$ and $q_{2}$, in the group $c_{2}$
each person is asked $q_{2}$ and $q_{3}$, etc. In quantum mechanics
the matrix above could describe the well-known EPR/Bell paradigm with
two entangled spin-half particles: $q_{1}$ and $q_{3}$ correspond
to the two axes along which Alice measures spins in her particle,
while $q_{2}$ and $q_{4}$ correspond to the two axes analogously
used by Bob in his particle. 

Let us use the notation $R_{i}^{j}$ for the outcome of a measurement
of $q_{i}$ in context $c_{j}$: 
\begin{center}
\begin{tabular}{|c|c|c|c|c|}
\hline 
$R_{1}^{1}$ & $R_{2}^{1}$ &  &  & $c_{1}$\tabularnewline
\hline 
 & $R_{2}^{2}$ & $R_{3}^{2}$ &  & $c_{2}$\tabularnewline
\hline 
 &  & $R_{3}^{3}$ & $R_{4}^{3}$ & $c_{3}$\tabularnewline
\hline 
$R_{1}^{4}$ &  &  & $R_{4}^{4}$ & $c_{4}$\tabularnewline
\hline 
$q_{1}$ & $q_{2}$ & $q_{3}$ & $q_{4}$ & \multicolumn{1}{c}{}\tabularnewline
\cline{1-4} 
\end{tabular}$\:$.
\par\end{center}

\noindent Since the values of $R_{i}^{i}$ and $R_{i\oplus1}^{i}$
are empirically paired (two responses given by the same person, or
the measurements by Bob and Alice made simultaneously), the random
variables in each row of the matrix are jointly distributed. This
is not true for measurements made in different contexts: their joint
distribution is undefined, and we call them stochastically unrelated
(not to be confused with being stochastically independent, which is
a special case of being jointly distributed). In particular, $R_{i}^{j}$
and $R_{i}^{j'}$ measuring the same property $q_{i}$ in two different
contexts are stochastically unrelated. In the Kolmogorovian language,
$R_{i}^{j}$ and $R_{i}^{j'}$ are defined on different domain probability
spaces. It is therefore impossible to say that $R_{i}^{j}=R_{i}^{j'}$,
because equality is a special case of joint distribution. 

It is clear now that CPT imposes no constraints whatever on the row-wise
joint distributions. The CHSH/Fine inequality (\ref{eq: CHSH}) cannot
be derived for this matrix of contextually labelled random variables.
However, it can be derived as a solution for the following problem:
find necessary and sufficient conditions for the existence of a jointly
distributed quadruple of random variables $\left(S_{1},S_{2},S_{3},S_{4}\right)$
such that 
\begin{equation}
\left(S_{i},S_{i\oplus1}\right)\textnormal{ has the same distribution as }\left(R_{i}^{i},R_{i\oplus1}^{i}\right),
\end{equation}
for $i=1,\ldots,4$. Such a vector $\left(S_{1},S_{2},S_{3},S_{4}\right)$
is called a ``reduced coupling'' of the stochastically unrelated
pairs $\left(R_{i}^{i},R_{i\oplus1}^{i}\right)$ (Dzhafarov \& Kujala,
2016b). 

The reduced coupling $\left(S_{1},S_{2},S_{3},S_{4}\right)$ is merely
a shortcut for describing a special case of what we call a $\mathsf{C}$-coupling
(Dzhafarov, Cervantes, \& Kujala, 2017; Dzhafarov \& Kujala, 2016a,
2017a, 2017b). $\mathsf{C}$ is some property of a pair of random
variables, and a $\mathsf{C}$-coupling of the pairs $\left(R_{i}^{i},R_{i\oplus1}^{i}\right)$
in our example is a jointly distributed octuple of random variables
$\left(S_{1}^{1},S_{2}^{1},S_{2}^{2},S_{3}^{2},S_{3}^{3},S_{4}^{3},S_{4}^{4},S_{1}^{4}\right)$
such that, for $i=1,\ldots,4$,
\begin{equation}
\left(S_{i}^{i},S_{i\oplus1}^{i}\right)\textnormal{ has the same distribution as }\left(R_{i}^{i},R_{i\oplus1}^{i}\right),
\end{equation}
and, in addition,
\begin{equation}
\left(S_{i\oplus1}^{i},S_{i\oplus1}^{i\oplus1}\right)\textnormal{ satisfies property }\mathsf{C}.
\end{equation}
The reduced coupling is the one defined by $\mathsf{C}$ with the
meaning ``are equal with probability 1'' (applied to pairs of random
variables). More generally, in the Contextuality-by-Default theory,
$\mathsf{C}$ is chosen to mean ``are equal with maximal possible
probability.'' For this choice of $\mathsf{C}$, the criterion for
the existence of a $\mathsf{C}$-coupling is
\begin{equation}
\begin{array}{r}
\max_{j=1,\ldots,4}\left|\sum_{i=1}^{4}\left\langle R_{i}^{i}R_{i\oplus1}^{i}\right\rangle -2\left\langle R_{j}^{j}R_{j\oplus1}^{j}\right\rangle \right|\\
\\
\leq2+\sum_{i=1}^{4}\left|\left\langle R_{i\oplus1}^{i}\right\rangle -\left\langle R_{i\oplus1}^{i\oplus1}\right\rangle \right|,
\end{array}\label{eq: criterion for cyclic 4}
\end{equation}
a useful generalization of CHSH/Fine inequality (\ref{eq: CHSH})
(Dzhafarov \& Kujala, 2016a; Dzhafarov, Kujala, \& Larsson, 2015;
Kujala \& Dzhafarov, 2016; Kujala, Dzhafarov, \& Larsson, 2015). A
system of random variables for which a $\mathsf{C}$-coupling exists
(does not exist) is called $\mathsf{C}$-noncontextual (respectively,
$\mathsf{C}$-contextual).\footnote{To avoid terminological confusion, in the Contextuality-by-Default
theory each random variable is contextually labeled (``by default''),
but a system of contextually labelled random variables can be $\mathsf{C}$-noncontextual
or $\mathsf{C}$-contextual. }

That context is part of the identity of a random variable is the departure
point for the Contextuality-by-Default theory, the term ``identity''
being understood in the Kolmogorovian sense, as the measurable function
from a domain probability space to a codomain measurable space (for
detailed explanations, see Dzhafarov \& Kujala, 2016a, 2017a). One
advantage provided by this contextual identification is that it allows
for the possibility that random variables measuring the same property
in different contexts, such as $R_{2}^{1}$ and $R_{2}^{2}$ in our
example, are differently distributed. This can happen, e.g., if one
of the two questions posed to a person influences her response to
the other question, or if Alice can signal to Bob and thereby change
his recordings. With non-contextual labeling, such as $R_{1},R_{2},R_{3},R_{4}$
in the opening formulation, to express the same fact one would have
to say that $R_{2}$ is differently distributed depending on whether
``it'' is recorded together with $R_{1}$ or $R_{3}$. This is at
best an abuse of language, if not outright nonsensical, as the distribution
of $R_{2}$ is part of its identity. 

\section{\label{sec:On-Statement-2}On Statement 2}
\begin{quotation}
``In classical probability the joint occurrence of two events is
commutative, but we know from quantum mechanics and from behavioral
experiments that the order of two events generally matters for their
joint probability.''
\end{quotation}
Thus, we read in Trueblood and Busemeyer (2011) that ``the classical
probability model has difficulty accounting for order effects because
the commutative property holds'' (p. 1527). And in Wang and Busemeyer
(2015): ``Classical probability theory has difficulty explaining
order effects because events are represented as sets and are commutative,
so the joint probability of events A and B is the same for the order
of `A and B' and the order of `B and A'{}'' (p. 2). Quotes like these
are numerous, but it should be noted that Busemeyer and colleagues
carefully qualify their criticism of CPT. They acknowledge that models
based on CPT can be formulated for such empirical phenomena as order
effects, but their presentation implies that these CPT-based models
have to be contrived. According to these authors, the only way CPT
can handle order effects is by using the Ballentine (1986) type conditionalization:
order of events (``B follows A'' and ``A follows B'') is considered
a random event conditioning probabilities of responses. We too consider
this construction awkward (Dzhafarov \& Kujala, 2014b; see also Section
\ref{sec:On-Statement-3}), but it is not the only one within the
framework of CPT: the Contextuality-by-Default approach provides another
way, one that is both simple and universally applicable.

Let us precede our discussion by pointing out that CPT would indeed
be a singularly helpless exercise if it lacked natural ways to depict
the difference between an ordered pairs of observations $\left(a,b\right)$
and an unordered two-element set $\{a,b\}$. The difference between
the two is obvious on the basic set-theoretic level: an ordered pair
$\left(a,b\right)$ is an abbreviation for the set $\left\{ \left\{ a,1\right\} ,\left\{ b,2\right\} \right\} $,
or $\left\{ a,\left\{ a,b\right\} \right\} $, because of which $\left(a,b\right)$
and $\left(b,a\right)$ are different sets, unless $a=b$, and $\left(a,a\right)=\left\{ \left\{ a,1\right\} ,\left\{ a,2\right\} \right\} $
is different from $\left\{ a,a\right\} =\left\{ a\right\} $. Moreover,
since an ordered pair is merely a simple case of a process (indexed
set), the logic of Statement 2 implies that CPT should resort to contrived
constructions when dealing with random processes that are not exchangeable.
A statement from Bruza, Wang, and Busemeyer (2015) may help in recognizing
that the order effects are a non-issue for CPT. The statement is that
in CPT ``the intersection of events is always defined and events
always commute, even if the events are distinguished by time (e.g.,
\textquoteleft $A$ at time 1\textquoteright{} and \textquoteleft $B$
at time 2\textquoteright{} is equivalent to \textquoteleft $B$ at
time 2\textquoteright{} and \textquoteleft $A$ at time 1\textquoteright )''
(p. 387). For ordered events it is precisely this commutativity, of
\textquoteleft $A$ at time 1\textquoteright{} and \textquoteleft $B$
at time 2\textquoteright , that holds in CPT (and classical logic).
In this form it is unchallengeable and cannot lead to any problems.
An issue is created when one compares \textquoteleft $A$ at time
1\textquoteright{} and \textquoteleft $B$ at time 2\textquoteright{}
to \textquoteleft $B$ at time 1\textquoteright{} and \textquoteleft $A$
at time 2\textquoteright , two conjunctions that are two different
events that need not have the same probability in CPT.

However, the approach offered by the Contextuality-by-Default theory
does not consist in labeling events. Rather, it uses a more versatile
labeling of random variables (although the two are essentially equivalent
in the simple case of two consecutive events). To understand the logic
of the approach, consider the probabilistic identity of responses
$R_{q}$ to some question $q$. Its domain probability space can be
thought of as a set $X$ of potential responders to the question $q$,
with some probability measure $\mu$ imposed on its power set (treated
as sigma-algebra). Let the possible values of $R_{q}$ be Yes/No.
Its distribution then is defined by
\begin{equation}
\Pr\left[R_{q}=\textnormal{Yes}\right]=\mu\left(\left\{ x\in X:x\textnormal{ responds Yes to }q\right\} \right).
\end{equation}
By construction, $q$ is part of the identity of $R_{q}$, so if $q$
is replaced with another question $q'$, the random variable $R_{q}$
will be replaced with another random variable $R_{q'}$. Probability
theory allows this new random variable to have another distribution,
but, of course, being an abstract mathematical theory, it does not
predict what the distributions of $R_{q}$ and $R_{q'}$ can be: such
a prediction is up to an empirical theory dealing with people's substantive
knowledge of questions and answers. 

Consider now two questions that have identical formulation but are
asked in different tones of voice or with different noise or images
in the background; or two questions that have the same content but
differ in how they are formulated (e.g., ``Is it 11 am now?'' versus
``Is 11 am the correct time at this moment?''). The usual experimental
design, if one is interested in such differences, would be to partition
$X$ into two subsets $X^{1}$ and $X^{2}$, asking the question $q$
in one form of the members of $X^{1}$ and in another form of the
members of $X^{2}$. From the point of view of abstract probability
theory, whatever the difference between the two questions substantively,
formally the responses to them are two different random variables
defined on two different domain probability spaces. They are, therefore,
stochastically unrelated. One can choose (based on one's substantive,
non-mathematical understanding of questions and answers) to consider
the differences in formulations or in the tone of voice to be part
of the questions themselves (in which case one will deal with random
variables denoted $R_{q}$ and $R_{q'}$) or to formalize the differences
as different contexts in which one and the same question is asked
(in which case one will present the random variables as $R_{q}^{c}$
and $R_{q}^{c'}$). The two representations are interchangeable, but
the latter one is preferable because, in accordance with the principles
of the Contextuality-by-Default theory, it encodes the stochastic
unrelatedness of $R_{q}^{c}$ and $R_{q}^{c'}$ in the very notation.\footnote{A combined notation, such as $R_{q}^{c}$ and $R_{q'}^{c'}$, is possible
too, as discussed by Dzhafarov and Kujala (2015), but it is less interesting
for subsequent contextuality analysis. As a general principle, it
is possible but counterproductive to include contexts as part of contents
of random variables: strictly separating the two is essential for
any contextuality analysis.} 

Using different orders of two questions has precisely the same logical
status as differences in the tone of voice or background noise: it
creates two pairs of jointly distributed random variables that are
stochastically unrelated to each other. The set $X$ is partitioned
into two subsets $X^{AB}$ and $X^{BA}$, corresponding to the two
orders, $\left(q_{A},q_{B}\right)$ and $\left(q_{B},q_{A}\right)$.
The random variables defined on these subsets and corresponding to
a given question, say $q_{A}$, can have different distributions.
The latter is exceedingly obvious if one uses specially chosen questions.
Consider, e.g., $q_{A}=$``Is this the first question I am asking?''
and $q_{B}=$``Is this the second question I am asking?'', asked
in two different orders. 

By analogy with two forms of the same question, one can now proceed
in several different ways, but the one most informative for contextuality
analysis is as follows. We define a jointly distributed pair $\left(R_{A}^{AB},R_{B}^{AB}\right)$
with
\begin{equation}
\begin{array}{c}
\Pr\left[R_{A}^{AB}=\textnormal{Yes}\right]=\mu\left(X_{A}^{AB}\right),\\
\\
\Pr\left[R_{B}^{AB}=\textnormal{Yes}\right]=\mu\left(X_{B}^{AB}\right),\\
\\
\Pr\left[R_{A}^{AB}=\textnormal{Yes \& }R_{B}^{AB}=\textnormal{Yes}\right]=\mu\left(X_{A}^{AB}\cap X_{B}^{AB}\right),
\end{array}\label{eq: AB}
\end{equation}
where
\begin{equation}
\begin{array}{c}
X_{A}^{AB}=\left\{ x\in X^{AB}:x\textnormal{ responds Yes to }q_{A}\right\} ,\\
\\
X_{B}^{AB}=\left\{ x\in X^{AB}:x\textnormal{ responds Yes to }q_{B}\right\} .
\end{array}
\end{equation}
The joint distribution for $\left(R_{A}^{BA},R_{B}^{BA}\right)$,
stochastically unrelated to the previous pair, is defined similarly,
and can be arbitrarily different from (\ref{eq: AB}). 

Using the Contextuality-by-Default representation, the system of the
random variables just defined is
\begin{center}
\begin{tabular}{|c|c|c|}
\hline 
$R_{A}^{AB}$ & $R_{B}^{AB}$ & $c_{AB}=\left(q_{A},q_{B}\right)$\tabularnewline
\hline 
$R_{A}^{BA}$ & $R_{B}^{BA}$ & $c_{BA}=\left(q_{B},q_{A}\right)$\tabularnewline
\hline 
$q_{A}$ & $q_{B}$ & \multicolumn{1}{c}{}\tabularnewline
\cline{1-2} 
\end{tabular}$\:$.
\par\end{center}

\noindent We can now choose some property $\mathsf{C}$ for pairs
of random variables, as explained in the previous section, and ask
whether the system above has a $\mathsf{C}$-coupling (or, in the
terminology of Contextuality-by-Default, whether it is $\mathsf{C}$-noncontextual).
With $\mathsf{C}$ chosen to mean ``are equal with maximal possible
probability,'' such a $\mathsf{C}$-coupling exists if and only if
(Dzhafarov \& Kujala, 2016a; Dzhafarov, Zhang, \& Kujala, 2015)

\begin{equation}
\begin{array}{l}
\left|\left\langle R_{A}^{AB}R_{B}^{AB}\right\rangle -\left\langle R_{A}^{BA}R_{B}^{BA}\right\rangle \right|\\
\\
\leq\left|\left\langle R_{A}^{AB}\right\rangle -\left\langle R_{A}^{BA}\right\rangle \right|+\left|\left\langle R_{B}^{AB}\right\rangle -\left\langle R_{B}^{BA}\right\rangle \right|.
\end{array}\label{eq: cyclic 2}
\end{equation}
Here, Yes and No responses have been encoded as $+1$ and $-1$, respectively.
The remarkable QQ equality discovered by Wang and Busemeyer (2013)
is equivalent to saying that the left-hand side expression in (\ref{eq: cyclic 2})
is zero, from which it follows that according to this law this system
of random variables is noncontextual.\footnote{Wang and Busemeyer (2013) are right about CPT (they call it ``Bayesian'')
not being able to predict the QQ equality, although they seem to consider
this a deficiency rather than a hallmark of any abstract mathematical
theory. Any prediction derived in a special-purpose theory, such as
the quantum formalism used by Wang and Busemeyer to derive the QQ
equality, can be (and always is, eventually) fully expressed in the
language of CPT, making it possible to relate the prediction to empirical
data.} See Dzhafarov, Kujala, Cervantes, Zhang, and Jones (2016) and Dzhafarov,
Zhang, and Kujala (2015) for a detailed discussion. 

\section{\label{sec:On-Statement-3}On Statement 3}
\begin{quotation}
``Classical probability is additive (equivalently, obeys the law
of total probability), but we know from quantum mechanics and from
behavioral experiments that this additivity (the law of total probability)
can be violated.''
\end{quotation}
Additivity, expressed in the language of random variables, is that
if $A$ and $B$ are disjoint events in the codomain space of a random
variable $R$, then
\begin{equation}
\Pr\left[R\in A\cup B\right]=\Pr\left[R\in A\right]+\Pr\left[R\in B\right].\label{eq: additivity}
\end{equation}
This principle is sometimes analyzed in an equivalent form, referred
to as the ``law of total probability'': as a consequence of the
additivity above and the set-theoretic distributivity, for any $C$
in the codomain space of $R$,
\begin{equation}
\begin{array}{l}
\Pr\left[R\in C\cap\left(A\cup B\right)\right]\\
\\
=\Pr\left[R\in C\cap A\right]+\Pr\left[R\in C\cap B\right].
\end{array}\label{eq: additivity+distributivity}
\end{equation}
Set-theoretic distributivity being an integral part of CPT, and the
two equalities above understood as belonging to CPT, it is logically
impossible to claim that one of them can be violated without stating
the same for the other. This should be kept in mind when encountering
statements about violations of the ``\emph{classical law}'' of total
probability. To give examples: ``It was shown that FTP {[}the formula
of total probability{]} (and hence classical probability theory) is
violated in some experiments on recognition of ambiguous pictures''
(Khrennikov, 2010, p. 90), and ``One can find evidence of violation
of laws of classical probability theory, e.g., in violation of the
law of total probability'' (Khrennikov \& Basieva, 2014, p. 105). 

We will discuss here the basic form (\ref{eq: additivity}) only. 

The claim of violations of this law in quantum mechanics comes from
the double-slit experiment. We consider it in the following version:
a source of particles emits them into a barrier with two slits (left
and right, each of which can be closed or open), and a detector of
the particles occupies a small area behind this barrier. One considers
the probability with which an emitted particle reaches the detector,
and discovers that this probability, when both slits are open, is
not equal to (depending on the detector's location, can be greater
or smaller than) the sum of these probabilities recorded with only
the left slit open and with only the right slit open. Richard Feynman
is often quoted as saying that this is ``a phenomenon which is impossible,
\emph{absolutely} impossible, to explain in any classical way'' (Feynman,
Leighton, \& Sands, 1975, Section 37-1). The words ``in any classical
way'' in this quote are commonly interpreted as ``by means of CPT.''
This interpretation may be correct, but it is also possible that Feynman
meant that this phenomenon cannot be explained by means of classical
mechanics, and that he viewed quantum probabilities as a special-purpose
theory for computing probabilities in a specific physical situation.
The second quote from Feynman (1951) given at the end of our introductory
section seems to agree with this interpretation. 

Whatever the case with Feynman, Ballentine (1986) presents a systematic
analysis of the double-slit experiment in terms of CPT, and argues
that the two are perfectly compatible if one treats the probabilities
in (\ref{eq: additivity}) as conditional ones, conditioned on three
different events. Translating this into the language of random variables,
Ballentine's solution is to rewrite (\ref{eq: additivity}) as 
\begin{equation}
\begin{array}{l}
\Pr\left[R\in A\cup B\,|\,Q=c_{\circ\circ}\right]\\
\\
\overset{generally}{\not=}\Pr\left[R\in A\,|\,Q=c_{\circ\times}\right]+\Pr\left[R\in B\,|\,Q=c_{\times\circ}\right],
\end{array}\label{eq: Ballentine}
\end{equation}
where $Q$ is a random variable indicating which of the slits is open
and which is closed: $c_{\circ\circ}$ means that both are open, $c_{\circ\times}$
means that only the left one is open, and $c_{\times\circ}$ means
the opposite. Clearly, the three conditioning values $c_{\circ\circ},c_{\circ\times},c_{\times\circ}$
are distinct and mutually exclusive (as with any distinct values of
any random variable), whence no equality in (\ref{eq: Ballentine})
should generally be expected.\footnote{In a personal communication (April 2018), Jerome Busemeyer explained
to us that when he and his colleagues speak of violations of the total
probability law they do not mean that (\ref{eq: additivity}) or (\ref{eq: additivity+distributivity})
fail to hold. Rather they mean that the probability of an event can
be different depending on what other events it is recorded together
with. This can be interpreted as a position very close to Contextuality-by-Default,
or to Ballentine's conceptualization (\ref{eq: Ballentine}), or even
as the possibility that, say, $C$ when one measures the probability
of $R\in C\cap A$ may be a different event from $C$ when one measures
the probability of $R\in C\cap\left(A\cup B\right)$ (the latter case
being formulated as measuring the probability of ``$C$ alone''
if $B=\textnormal{not}A$). While welcoming this clarification, one
should note that its implication is that the law of total probability
as a formula of CPT \emph{is not violated}. Ballentine's inequality
(\ref{eq: Ballentine}) is a correct CPT formula, and something like
$\Pr\left[R\in C_{1}\cap\left(A\cup B\right)\right]=\Pr\left[R\in C_{2}\cap A\right]+\Pr\left[R\in C_{3}\cap B\right]$
is not a correct CPT formula. (The latter, incidentally, is a good
example of why labeling of events, rather than of random variables,
is not a good solution: shall one, in addition to $C$, also differently
label $A$ and $B$ on the right and on the left of the formula?)} 

Dzhafarov and Kujala (2014b) call this approach ``conditionalization,''
in relation to Avis, Fischer, Hilbert, and Khrennikov (2009) where
it was used systematically (see also Khrennikov, 2006, 2015b). It
is true that if $R$ and $Q$ are jointly distributed, then $R$ conditioned
on some value of $Q$ and $R$ conditioned on another value of $Q$
are two random variables that possess no joint distributions, i.e.,
are stochastically unrelated. This means that conditionalization is
a special case of contextual labeling, in fact an instructive case
for introducing the notion of stochastic unrelatedness (Dzhafarov
\& Kujala, 2014b, 2016b). However, the choice among conditions $c_{\circ\circ},c_{\circ\times},c_{\times\circ}$
need not be random. One can conduct an experiment with both slits
open for a year, then for another year with the left slit closed,
and so on. This should not change anything in the analysis of the
double-slit experiment. As we mentioned in the introductory section,
Khrennikov (2009a) pointed out the difference between contextual and
conditional probabilities in presenting his general V\"axj\"o contextual
model.\footnote{In some of his later work, however, Khrennikov seems to have abandoned
this distinction and adopted a rigorous version of Ballentine's view
(Avis, Fischer, Hilbert, \& Khrennikov, 2009; Khrennikov, 2015a, 2015b;
see Dzhafarov \& Kujala, 2014b, for a critical discussion of this
position).}

The analysis of the double-slit experiment within the framework of
the Contextuality-by-Default theory begins with identifying the random
variables in play, their contexts and the properties they measure.
The contexts are the same as in Ballentine's analysis, $c_{\circ\circ},c_{\circ\times},c_{\times\circ}$,
but for completeness we add $c_{\times\times}$ (both slits closed).
The measured properties are identified by which slit one considers
(left or right) and by its condition (open or closed): $q_{\circ\cdot}$
(left slit, open), $q_{\cdot\circ}$(right slit, open), and similarly
for the closed slits, $q_{\times\cdot},q_{\cdot\times}$. This creates
eight random variables that we can arrange as follows: 
\begin{center}
\begin{tabular}{|c|c|c|c|c|}
\hline 
$R_{\circ\cdot}^{\circ\circ}$ & $R_{\cdot\circ}^{\circ\circ}$ &  &  & $c_{\circ\circ}$\tabularnewline
\hline 
 & $R_{\cdot\circ}^{\times\circ}$ & $R_{\times\cdot}^{\times\circ}$ &  & $c_{\times\circ}$\tabularnewline
\hline 
 &  & $R_{\times\cdot}^{\times\times}$ & $R_{\cdot\times}^{\times\times}$ & $c_{\times\times}$\tabularnewline
\hline 
$R_{\circ\cdot}^{\circ\times}$ &  &  & $R_{\cdot\times}^{\circ\times}$ & $c_{\circ\times}$\tabularnewline
\hline 
$q_{\circ\cdot}$ & $q_{\cdot\circ}$ & $q_{\times\cdot}$ & $q_{\cdot\times}$ & \multicolumn{1}{c}{}\tabularnewline
\cline{1-4} 
\end{tabular}$\:$,
\par\end{center}

\noindent where $R_{q}^{c}$ can be interpreted as answering the
question ``Has the particle just emitted reached the detector through
the slit $q$ under condition $c$?'' The possible values of $R_{q}^{c}$
are Yes and No. Thus, $R_{\times\cdot}^{\times\circ}=$Yes means that
a particle reached the detector having passed through the closed slit
on the left when the right slit is open. Physics (not probability
theory) tells us that the probability of this happening is zero. 

Our arrangement of the random variables shows that, surprisingly,
the system they comprise is formally a cyclic system of rank 4 (Dzhafarov
\& Kujala, 2016a; Dzhafarov, Kujala, \& Larsson, 2015; Kujala, Dzhafarov,
\& Larsson, 2015). It is the same system as the one in the simplest
EPR/Bell ``Alice-Bob'' paradigm, described in our analysis of Statement
1. If one uses one's knowledge that no particle can reach the detector
through a closed slit, then the joint distributions of all context-sharing
pairs of random variables (the rows of the matrix above) are defined
by the following joint and marginal probabilities:
\begin{center}
\begin{tabular}{c|c|c}
$c_{\circ\times}$ & $R_{\cdot\times}^{\circ\times}=$Yes & \tabularnewline
\hline 
$R_{\circ\cdot}^{\circ\times}=$Yes & $\ensuremath{0}$  & \multicolumn{1}{c|}{$p$}\tabularnewline
\hline 
 & $0$ & \tabularnewline
\cline{2-2} 
\end{tabular}$\:$,$\:$%
\begin{tabular}{c|c|c}
$c_{\times\times}$ & $R_{\cdot\times}^{\times\times}=$Yes & \tabularnewline
\hline 
$R_{\times\cdot}^{\times\times}=$Yes & $\ensuremath{0}$  & \multicolumn{1}{c|}{$0$}\tabularnewline
\hline 
 & $0$ & \tabularnewline
\cline{2-2} 
\end{tabular}$\:$,\smallskip{}
\par\end{center}

\begin{center}
\begin{tabular}{c|c|c}
$c_{\times\circ}$ & $R_{\cdot\circ}^{\times\circ}=$Yes & \tabularnewline
\hline 
$R_{\times\cdot}^{\times\circ}=$Yes & $\ensuremath{0}$  & \multicolumn{1}{c|}{$0$}\tabularnewline
\hline 
 & $q$ & \tabularnewline
\cline{2-2} 
\end{tabular} $\:$,$\:$%
\begin{tabular}{c|c|c}
$c_{\circ\circ}$ & $R_{\cdot\circ}^{\circ\circ}=$Yes & \tabularnewline
\hline 
$R_{\circ\cdot}^{\circ\circ}=$Yes & $\ensuremath{r'}$  & \multicolumn{1}{c|}{$r'+p'$}\tabularnewline
\hline 
 & $r'+q'$ & \tabularnewline
\cline{2-2} 
\end{tabular}$\:$,
\par\end{center}

\begin{center}
\smallskip{}
\par\end{center}

\noindent where $p,q,p',q',r'$ are some probabilities. The physical
interpretation of the joint distribution for $c_{\circ\circ}$ compared
to that for, say, $c_{\circ\times}$ is that, somehow, the way particles
reach the detector having passed through the open left slit may be
different depending on whether the right slit is open or closed. A
physicist may tell us that this is because of the particle-wave duality
and wave interference, but this is irrelevant for the probabilistic
analysis. 

It is interesting to see whether the system just described is $\mathsf{C}$-noncontextual
(has a $\mathsf{C}$-coupling) with $\mathsf{C}$=``are equal with
maximal possible probability.'' The application of the general criterion
(\ref{eq: criterion for cyclic 4}) for noncontextuality of such a
system yields
\begin{equation}
\begin{array}{r}
\left(\left(1-2p\right)+\left(1\right)+\left(1-2q\right)+(1-2p'-2q')\right)\\
\\
-2\min\left(\left(1-2p\right),\left(1\right),\left(1-2q\right),(1-2p'-2q')\right)\\
\\
\leq2+2|p-p'-r'|+|-1+1|\\
\\
+2|q-q'-r'|+|-1+1|,
\end{array}
\end{equation}
where we have assumed that the detector is so small that the probabilities
$1-2p,1-2q,1-2p'-2q'$ are all positive. By simple algebra one can
show that this inequality is always satisfied, that is, the double-slit
system is $\mathsf{C}$-noncontextual. 

Comparing again $c_{\circ\circ}$ with $c_{\circ\times}$ (or $c_{\times\circ}$),
the noncontextuality just established means, within the framework
of Contextuality-by-Default, that the influence exerted by the state
of a slit (open or closed) upon how the particles reach the detector
having passed through the other slit is of a ``direct cross-influence''
nature, with no contextuality proper (Dzhafarov, 2017; Dzhafarov,
Cervantes, \& Kujala, 2017; Dzhafarov \& Kujala, 2017a). For a detailed
contextuality analysis of the double-slit experiment see Dzhafarov
and Kujala (2018), where it is also shown that a system with more
than two slits may very well be contextual.

\section{\label{sec:Concluding-remarks}Concluding remarks}

In this concluding section we will briefly address four commonly raised
concerns about the contextual notation and the principle that random
variables in different contexts are different (and stochastically
unrelated).\bigskip{}

\paragraph{Question: }

In empirical situations where the contexts are known not to influence
a measurement directly (like in the EPR/Bell Alice-Bob paradigm with
spacelike separation of the measurements), what ``causes'' the random
variable representing this measurement to change its identity?

\paragraph{Answer: }

The identity of a random variable is determined by its own distribution
and also by the joint distribution of this random variable with all
other random variables in the same context. Therefore, any change
in these other variables ``automatically'' changes its identity.
Here is a simple analogy. A person $P$ is in a room with other people.
$P$ has some characteristics, such as ``she is kind,'' or ``she
is tall.'' It is possible that she is the tallest person in the room,
in which case she is also characterized by this fact. The statement
``she is the tallest person in the room'' therefore describes a
property of $P$, part of her identity in addition to her being kind
and tall. If one of the other people leaves, and someone enters who
is taller than $P$, she ``automatically'' changes her identity,
as she ceases to be the tallest person in the room. This change in
$P$ occurs even if she is not aware of the change in the room, or
the room is so large that there are no physical means for her to notice
this. As acknowledged in Dzhafarov and Kujala (2014a), an application
of this general argument to quantum phenomena may be regarded as paralleling
an argument made by Bohr (1935) in reply to the EPR paper (Einstein,
Podolsky, \& Rosen, 1935).\bigskip{}

\paragraph{Question:}

If every condition recorded together with a random variable can be
considered part of its context, does this not mean that any two realizations
of the same random variable are in fact two different random variables,
stochastically unrelated to each other?

\paragraph{Answer: }

If the realizations are separately indexed, e.g., by the ordinal position
in a sequence of trials, each of them indeed must be viewed as a single
realization of a unique random variable. There is, however, a choice
of the point of view for subsequent analysis. One can view these unique
random variables as ones with different measured properties (trial
numbers) within a single context (sequence of trial numbers). Conversely,
one can view them as random variables measuring the same thing (e.g.,
they all measure the response of a person to a flash) but in different
contexts (trial number). We implicitly adopt the second point of view
when we speak of the sequence as one of different realizations of
``the same'' random variable. See Dzhafarov and Kujala (2015, 2016a)
for detailed discussions. In the case considered, the choice of one
of the two points of view makes no difference. If the realizations
are treated as context-sharing, the random variables are jointly distributed,
but this joint distribution is manifested in a single realization
only. We need additional assumptions to reconstruct it, such as stochastic
independence, ergodicity, martingale property, etc. If the realizations
are treated as measuring the same property in different contexts,
they are pairwise stochastically unrelated, and we need to couple
them. The choice of a coupling here amounts to adopting the same additional
assumptions.\bigskip{}

\paragraph{Question: }

Is the contextual labeling with stochastic unrelatedness really classical,
in the Kolmogorovian sense? 

\paragraph{Answer: }

It is a matter of definition and understanding of history. In some
publications one of us and Janne Kujala called our approach a ``qualified''
Kolmogorovian theory (Dzhafarov \& Kujala, 2014c), and it can be presented
in a way that sets it aside from a standard account of CPT (as in
Dzhafarov \& Kujala, 2016a). However, we prefer to speak of Contextuality-by-Default
as part of the Kolmogorovian probability theory, with a greater emphasis
on multiple freely introducible domain probability spaces, stochastically
unrelated random variables defined on these spaces, and their couplings
understood as placing their copies on the same domain space. This
preference is based on our disbelief that Kolmogorov himself and the
brilliant probabilists working in his language could have overlooked
the obvious fact that there cannot exist a joint distribution of all
imaginable random variables. In his celebrated little book, Kolmogorov
(1956, §2 of Chapter 1) discusses empirical applications of his mathematical
theory, and in doing so confines his consideration to a single experiment
(corresponding, in our language, to a single context). He may have
erroneously thought this was the only realistic or interesting application.
This reading of Kolmogorov is also advocated in Khrennikov (2009b).\bigskip{}

\paragraph{Question: }

If, however, one posits that in any application of CPT all random
variables involved are defined on a single domain probability space,
would not then the claim of the inadequacy of CPT \emph{thus understood}
be justified?

\paragraph{Answer: }

The issue of the existence of a joint distribution for all random
variables involved in a given application is not as critical as the
issue of \emph{what} random variables are involved. Statements 1,
2, and 3 considered above are based first and foremost on misidentifying
the random variables in play. Thus, the correct system of random variables
representing the question order experiment is
\begin{center}
$\mathsf{R}=$ %
\begin{tabular}{|c|c|c|}
\hline 
$R_{A}^{AB}$ & $R_{B}^{AB}$ & $c_{AB}=\left(q_{A},q_{B}\right)$\tabularnewline
\hline 
$R_{A}^{BA}$ & $R_{B}^{BA}$ & $c_{BA}=\left(q_{B},q_{A}\right)$\tabularnewline
\hline 
$q_{A}$ & $q_{B}$ & \multicolumn{1}{c}{}\tabularnewline
\cline{1-2} 
\end{tabular}.
\par\end{center}

\noindent It is simply unjustifiable to posit a priori that it can
be replaced with
\begin{center}
$\mathsf{R'}=$ %
\begin{tabular}{|c|c|c|}
\hline 
$R_{A}$ & $R_{B}$ & $c_{AB}=\left(q_{A},q_{B}\right)$\tabularnewline
\hline 
$R_{A}$ & $R_{B}$ & $c_{BA}=\left(q_{B},q_{A}\right)$\tabularnewline
\hline 
$q_{A}$ & $q_{B}$ & \multicolumn{1}{c}{}\tabularnewline
\cline{1-2} 
\end{tabular},
\par\end{center}

\noindent a system in which random variables do not change with context:
even if one ignores the logic of Contextuality-by-Default, there is
no rationale for assuming that contexts are irrelevant, because in
this particular example one even knows that the distributions of $R_{A}^{AB}$
and $R_{A}^{BA}$ are different (which is the very ``question order
effect'' that makes this paradigm interesting). The situation here
is no different from someone deciding to replace $\mathsf{R}$ with 
\begin{center}
$\mathsf{R''}=$ %
\begin{tabular}{|c|c|c|}
\hline 
$R^{AB}$ & $R^{AB}$ & $c_{AB}=\left(q_{A},q_{B}\right)$\tabularnewline
\hline 
$R^{BA}$ & $R^{BA}$ & $c_{BA}=\left(q_{B},q_{A}\right)$\tabularnewline
\hline 
$q_{A}$ & $q_{B}$ & \multicolumn{1}{c}{}\tabularnewline
\cline{1-2} 
\end{tabular},
\par\end{center}

\noindent a system in which random variables do not change with content.

This reasoning applies even if one views the four random variables
in the original system $\mathsf{R}$ as having an unknown (and unknowable)
joint distribution. This amounts to informally identifying the system
with one of its possible couplings, and the construction of a $\mathsf{C}$-coupling
then can be presented as determining if this ``true'' joint distribution
could possibly be satisfying $\mathsf{C}$. One can check that our
analysis of the question order effect in Section \ref{sec:On-Statement-2}
would hold with no serious modifications if one adopted this language
(and similarly for the systems considered in Sections \ref{sec:On-Statement-1}
and \ref{sec:On-Statement-3}). 

There are, of course, good reasons not to use this language, except
as an informal version of the rigorous language of the Contextuality-by
Default theory (perhaps for the sake of conceptual or notational simplicity).
The assumption that any two random variables are jointly distributed
is mathematically untenable. It is untenable because, due to the transitivity
of the relation of being defined on the same domain probability space,
it implies the erroneous notion that there is a joint distribution
for all imaginable random variables (for reasons why this notion is
wrong, see Dzhafarov \& Kujala, 2014a, 2014b, 2017a). Within the framework
of CPT stochastically unrelated random variables must exist, making
it unjustifiable to assume without critical examination that all random
variables in a given set have a joint distribution.

\subsection*{Acknowledgments}

This research has been supported by AFOSR grant FA9550-14-1-0318 (E.D.)
and Purdue University Lynn Fellowship (M.K.). The authors are grateful
to Victor H. Cervantes for valuable critical suggestions, and to Tim
Maudlin for discussing with us Feynman's position on classical probability
theory. 

\section*{REFERENCES}

\setlength{\parindent}{0cm}\everypar={\hangindent=15pt}

Accardi, L. (1982). Foundations of quantum probability. Rendiconti
del Seminario Matematico, 40, 249-270.

Aerts, D. (2009). Quantum particles as conceptual entities: A possible
explanatory framework for quantum theory. Foundations of Science,
14, 361\textendash 411.

Aerts, D. (2014). Quantum theory and human perception of the macro-world.
Frontiers in Psychology, 5, 1-19. 

Aerts, D., \& Sozzo, S. (2011). Quantum structure in cognition: Why
and how concepts are entangled. In D. Song, M. Melucci, I. Frommholz,
P. Zhang, L. Wang, \& S. Arafat. (Eds.), Proceedings of the Quantum
Interaction Conference (pp. 116\textendash 127). Berlin: Springer.

Avis, D., Fischer, P., Hilbert, A., \& Khrennikov, A. (2009). Single,
complete, probability spaces consistent with EPR-Bohm-Bell experimental
data. In A. Khrennikov (Ed.), Foundations of Probability and Physics-5,
AIP Conference Proceedings 750 (pp. 294-301). Melville, New York:
AIP.

Ballentine, L.E. (1986). Probability theory in quantum mechanics.
American Journal of Physics, 54, 883-889.

Bell, J. (1964). On the Einstein-Podolsky-Rosen paradox. Physics,
1, 195-200. 

Broekaert, J., Basieva, I., Blasiak, P., \& Pothos, E.M. (2017). Quantum-like
dynamics applied to cognition: A consideration of available options.
Philosophical Transactions of the Royal Society A, 375: 20160387.

Bohr, N. (1935). Can quantum-mechanical description of physical reality
be considered complete? Physical Review, 48, 696-702

Bruza, P., Kitto, K., Nelson, D., \& McEvoy, C. (2009). Is there something
quantum-like about the human mental lexicon? Journal of Mathematical
Psychology, 53, 362\textendash 377.

Bruza, P., Wang, Z., \& Busemeyer, J.R. (2015). Quantum cognition:
A new theoretical approach to psychology. Trends in Cognitive Sciences
19, 383-393. 

Busemeyer, J.R., \& Bruza, P.D. (2012). Quantum Models of Cognition
and Decision. Cambridge: Cambridge University Press. 

Clauser, J.F., Horne, M.A., Shimony, A., \& Holt, R.A. (1969). Proposed
experiment to test local hidden-variable theories. Physical Review
Letters, 23, 880\textendash 884.

Costantini, D. (1993). A statistical analysis of the two-slit experiment:
Or some remarks on quantum probability. International Journal of Theoretical
Physics, 32, 2349\textendash 2362.

Dzhafarov, E.N. (2017). Replacing nothing with something special:
Contextuality-by-Default and dummy measurements. In A. Khrennikov
\& T. Bourama (Eds.), Quantum Foundations, Probability and Information.
Springer.

Dzhafarov, E.N., Cervantes, V.H., \& Kujala, J.V. (2017). Contextuality
in canonical systems of random variables. Philosophical Transactions
of the Royal Society A, 375: 20160389.

Dzhafarov, E.N., \& Kujala, J.V. (2014a). Contextuality is about identity
of random variables. Physica Scripta T163, 014009.

Dzhafarov, E.N., \& Kujala, J.V. (2014b). Embedding quantum into classical:
Contextualization vs conditionalization. PLoS One 9(3): e92818. doi:10.1371/journal.pone.0092818.

Dzhafarov, E.N., \& Kujala, J.V. (2014c). A qualified Kolmogorovian
account of probabilistic contextuality. In H. Atmanspacher, E. Haven,
K. Kitto, \& D. Raine (Eds.), Lecture Notes in Computer Science, 8369,
201-212.

Dzhafarov, E.N., \& Kujala, J.V. (2015). Conversations on contextuality.
In E.N. Dzhafarov, J.S. Jordan, R. Zhang, \& V.H Cervantes (Eds.),
Contextuality from Quantum Physics to Psychology (pp. 1-22). New Jersey:
World Scientific.

Dzhafarov, E.N., \& Kujala, J.V. (2016a). Context-content systems
of random variables: The contextuality-by-default theory. Journal
of Mathematical Psychology, 74, 11-33.

Dzhafarov, E.N., \& Kujala, J.V. (2016b). Probability, random variables,
and selectivity. In W. Batchelder, H. Colonius, E.N. Dzhafarov, \&
J. Myung (Eds.), The New Handbook of Mathematical Psychology (pp.
85-150). Cambridge: Cambridge University Press. 

Dzhafarov, E.N., \& Kujala, J.V. (2017a). Probabilistic foundations
of contextuality. Fortschritte der Physik - Progress of Physics, 65,
1-11.

Dzhafarov, E.N., \& Kujala, J.V. (2017b). Contextuality-by-Default
2.0: Systems with binary random variables. In J.A. de Barros, B. Coecke,
\& E. Pothos (Eds.), Lecture Notes in Computater Science, 10106, 16-32.

Dzhafarov, E.N., \& Kujala, J.V. (2018). Contextuality analysis of
the double slit experiment (with a glimpse into three slits). Entropy
20, 278; doi:10.3390/e20040278.

Dzhafarov, E.N., Kujala, J.V., \& Cervantes, V.H. (2016). Contextuality-by-Default:
A brief overview of ideas, concepts, and terminology. In H. Atmanspacher,
T. Filk, \& E. Pothos (Eds.), Lecture Notes in Computer Science, 9535,
12-23.

Dzhafarov, E.N., Kujala, J.V., Cervantes, V.H., Zhang, R., \& Jones,
M. (2016). On contextuality in behavioral data. Philosophical Transactions
of the Royal Society A, 374: 20150234.

Dzhafarov, E.N., Kujala, J.V., \& Larsson, J-Å. (2015) Contextuality
in three types of quantum-mechanical systems. Foundations of Physics,
7, 762-782.

Dzhafarov, E.N., Zhang, R., \& Kujala, J.V. (2015). Is there contextuality
in behavioral and social systems? Philosophical Transactions of the
Royal Society A, 374: 20150099.

Einstein, A., Podolsky, B., \& Rosen, N. (1935). Can quantum-mechanical
description of physical reality be considered complete? Physical Review,
47, 777-780.

Feynman, R.P. (1951). The concept of probability in quantum mechanics.
In J. Neyman (Ed.), Proceedings of the Second Berkeley Symposium on
Mathematical Statistics and Probability (pp. 533-541). Berkeley: University
of California Press.

Feynman, R.P., Leighton, R., \& Sands, M. (1975). The Feynman Lectures
on Physics. Reading, MA: Addison Wesley.

Filipp, S., \& Svozil, K. (2005). Tracing the bounds on Bell-type
inequalities. AIP Conference Proceedings, 750, 87\textendash 94. 

Fine, A. (1982). Hidden variables, joint probability, and the Bell
inequalities. Physical Review Letters, 48, 291\textendash 295.

Khrennikov, A. (2006). A formula of total probability with interference
term and the Hilbert space representation of the contextual Kolmogorovian
model. Theory of Probability and Applications, 51, 427\textendash 441.

Khrennikov, A. (2009a). Contextual Approach to Quantum Formalism.
Dordrecht: Springer.

Khrennikov, A. (2009b). Bell\textquoteright s inequality: Physics
meets probability. Information Science, 179, 492-504.

Khrennikov, A.Y. (2010). Ubiquitous Quantum Structure: From Psychology
to Finance. New York: Springer.

Khrennikov, A. (2015a). CHSH inequality: Quantum probabilities as
classical conditional probabilities. Foundations of Physics, 45, 711-725.

Khrennikov, A. (2015b). Two-slit experiment: Quantum and classical
probabilities. Physica Scripta, 90, 1-9.

Khrennikov, A. \& Basieva, I. (2014). Quantum(-like) decision making:
On validity of the Aumann theorem. In: H. Atmanspacher, C. Bergomi,
T. Filk, K. Kitto (Eds.) Lecture Notes in Computer Science, 8951,
105-118.

Kolmogorov, A.N. (1956). Foundations of the Theory of Probability.
New York: Chelsea Publishing Company.

Kujala, J.V., \& Dzhafarov, E.N. (2016). Proof of a conjecture on
contextuality in cyclic systems with binary variables. Foundations
of Physics, 46, 282-299.

Kujala, J.V., Dzhafarov, E.N., \& Larsson, J-Å. (2015). Necessary
and sufficient conditions for maximal noncontextuality in a broad
class of quantum mechanical systems. Physical Review Letters, 115,
150401.

Maudlin, T. (2013). Statistical Modeling, Causal Inference, and Social
Science. (2013, September 25). Classical probability does not apply
to quantum systems (causal inference edition) {[}Blog post{]}. Retrieved
from http://andrewgelman.com/2013/09/25/classical-probability-does-not-apply-to-quantum-systems-causal-inference-edition/.

Moreira, C., \& Wichert, A. (2016). Quantum probabilistic models revisited:
The case of disjunction effects in cognition. Frontiers in Physics,
4:26. doi: 10.3389/fphy.2016.00026.

Pothos, E.M., \& Busemeyer, J.R. (2013). Can quantum probability provide
a new direction for cognitive modeling? Behavioral and Brain Sciences,
36, 255\textendash 274. 

Trueblood, J.S., \& Busemeyer, J.R. (2011). A quantum probability
account of order effects in inference. Cognitive Science, 35, 1518\textendash 1552. 

Wang, Z., \& Busemeyer, J.R. (2013). A quantum question order model
supported by empirical tests of an a priori and precise prediction.
Topics in Cognitive Science, 5, 689\textendash 710.

Wang, Z., \& Busemeyer, J.R. (2015). Reintroducing the concept of
complementarity into psychology. Frontiers in Psychology, doi: 10.3389/fpsyg.2015.01822.

Yearsley, J.M., \& Pothos, E.M. (2014). Challenging the classical
notion of time in cognition: A quantum perspective. Proceedings of
the Royal Society B, 281: 20133056.
\end{document}